\documentclass[oneside,leqno]{article}
\usepackage{amssymb,amsmath,latexsym}
\def\gg{\mathfrak{g}}
\def\gh{\mathfrak{h}}

\def\gl{\mathfrak{l}}

\def\gn{\mathfrak{n}}

\def\gp{\mathfrak{p}}

\def\gs{\mathfrak{s}}

\def\C{\mathbb{C}}

\def\R{\mathbb{R}}

\def\Z{\mathbb{Z}}

\def\cT{\mathcal{T}}

\def\nek{\text{\hbox{$\simeq$ \kern-.95em \hbox{$/$ \kern.05em}}}}
\def\vep{\varepsilon}

\DeclareMathOperator{\Ker}{Ker}
\DeclareMathOperator{\Span}{Span}

\DeclareMathOperator{\Img}{Im}

\DeclareMathOperator{\supp}{supp}

\DeclareMathOperator{\tr}{tr}

\DeclareMathOperator{\hgt}{ht}

\def\cplus{\hbox{$\subset${\raise1.05pt\hbox{\kern -0.55em
${\scriptscriptstyle +}$}}\ }}
\def\bcplus{\hbox{$\supset${\raise1.05pt\hbox{\kern -0.55em
${\scriptscriptstyle +}$}}\ }}

\def\ctimes{\hbox{$\times${\raise1.1pt\hbox{\kern -0.27em
${\scriptscriptstyle |}$}}\ }}

\def\udarrow{\hbox{$\nearrow${\kern -0.97em$\searrow$}\ }}

\def\bctimes{\hbox{$\times${\raise1.1pt\hbox{\kern -.74em
${\scriptscriptstyle |}$}}\ }\,\,}

\newtheorem{theorem}[equation]{Theorem}

\newtheorem{proposition}[equation]{Proposition}

\newtheorem{definition}[equation]{Definition}

\title{Parabolic sets of roots}
\author{Ivan Dimitrov\thanks{\, Research partially supported by NSERC Discovery Grant and by FAPESP (2007/03735-0), Brazil} ,
Vyacheslav Futorny\thanks{\, Research partially supported by CNPq
(301743/2007-0) and by Fapesp (2005/60337-2), Brazil} , and
Dimitar Grantcharov}

\date{\empty}

\begin{document}

\maketitle

\begin{abstract}
We compare two combinatorial definitions of parabolic sets of roots. We show that these
definitions are equivalent for simple finite dimensional Lie algebras, affine Lie algebras, and toroidal
Lie algebras. In contrast, these definitions are not always equivalent for simple finite dimensional Lie superalgebras.

Key words (2000 MSC): Primary 17B20, Secondary 17B65, 17B67

\end{abstract}

\section*{Introduction}

Let $\gg$ be a Lie algebra or superalgebra over $\C$ and let $\gh$ be an abelian subalgebra of $\gg$ which acts semisimply on $\gg$. The decomposition
\begin{equation} \label{eq11}
\gg = \gg^0 \oplus \left( \bigoplus_{\alpha \in \Delta} \, \gg^\alpha \right),
\end{equation}
where $\Delta \subset \gh^* \setminus \{ 0 \}$ and $\gg^\alpha = \{ x \in \gg \, | \, [h, x] = \alpha(h) x {\text { for every }} h \in \gh\}$ is called root decomposition of $\gg$ with respect to $\gh$.
The set $\Delta$ is called the root system of $\gg$ with respect to $\gh$. Note that a priori $\Delta$ depends both on $\gg$ and $\gh$ but since we will always consider a
fixed pair $(\gg, \gh)$, we are not going to emphasize this dependence. A $\gg$--module $M$ is called a {\it weight module} if it is semisimple as an $\gh$--module, i.e.
$$
M = \bigoplus_{\mu \in \gh^*} M^\mu,  {\text { where }} M^\mu = \{ m \in M \, | \, h \cdot m = \mu(h) m {\text { for every }} h \in \gh\}.
$$
The elements $\mu \in \gh^*$ with $M^\mu \neq 0$ are called {\it weights} of $M$ and the corresponding  $M^\mu$ are called {\it weight spaces} of $M$. The set $\supp M$ of all weights of $M$ is called the {\it support} of $M$.

One of the major achievements of the study of representations of
finite dimensional reductive Lie algebras is the classification of
irreducible weight modules with finite dimensional weight spaces.
The first step in this classification is the Fernando-Futorny theorem,
proved by V. Futorny for classical Lie algebras, \cite{Fu5}, and by
S. Fernando for arbitrary simple Lie algebras, \cite{Fe}. It
states that every irredicible weight $\gg$--module $M$ with finite dimensional weight spaces
 is either torsion free
or there exist a proper parabolic subalgebra $\gp$ of $\gg$ and an
irreducible weight module $N$ with finite dimensional weight
spaces over the Levi component of $\gp$ such that $M$ is
isomorphic to the unique irreducible quotient of the $\gg$--module
parabolically induced from $N$.

Here are some  details about the proof.

Fix an irreducible weight $\gg$--module $M$ with finite dimensional weight spaces.
Let $P$ be the set of roots $\alpha$ for which either $\gg^\alpha$ acts locally nilpotently on $M$ or both $\gg^\alpha$ and $\gg^{-\alpha}$ act injectively on $M$.
The definition of $P$ implies that
\begin{equation} \label{P1}
 \Delta = P \cup -P.
\end{equation}
A theorem proved independently by Fernando, \cite{Fe}, and Kac, \cite{K2}, implies that $P$ is a closed subset of $\Delta$, i.e. that
\begin{equation} \label{P2}
\alpha, \beta \in P {\text { with }} \alpha + \beta \in \Delta {\text { implies }} \alpha + \beta \in P.
\end{equation}
Properties (\ref{P1}) and (\ref{P2}) imply that
\begin{equation} \label{parabolic}
\gp = \gh \oplus \left( \bigoplus_{\alpha \in P} \, \gg^\alpha \right)
\end{equation}
is a parabolic subalgebra of $\gg$. It equals $\gg$ if all root spaces $\gg^\alpha$ act locally nilpotently on $M$ or all root spaces $\gg^\alpha$
act injectively on $M$. In the former case $M$ is a finite dimensional module and in the latter case it is a torsion free module. Torsion free modules
were classified by Mathieu in \cite{M}. We now turn our attention to the case when $\gp$ is a proper subalgebra of $\gg$. If we denote by $\gn$ and $\gl$ the
nilpotent radical and the reductive part of $\gp$ respectively, then the vector space $N = M^\gn$ of $\gn$--invariants of $M$ is nonzero and has a structure
of both an $\gl$--module and a $\gp$--module. Moreover, $N$ is irreducible as a module over each of $\gl$ and $\gp$. Consider the parabolically induced
module
$$
M_\gp(N) = U(\gg) \otimes_{U(\gp)} N,
$$
where $U(\gg)$ and $U(\gp)$ denote the universal enveloping algebras of $\gg$ and $\gp$ respectively. The parabolic subalgebra $\gp$ induces a $\Z$--grading
\begin{equation} \label{grading}
\gg = \bigoplus_{i \in \Z} \gg_i \quad {\text{ such that }} \quad \gg_0 = \gl, {\text { and }} \gn = \bigoplus_{i > 0} \, \gg_i.
\end{equation}
 The grading (\ref{grading}) extends to a $\Z$--grading on $U(\gg)$ and, consequently,
on $M_\gp(N)$. Using the grading on $M_\gp(N)$ one proves easily that $M_\gp(N)$ has a unique maximal proper submodule and, hence, a unique irreducible quotient
$V_\gp(N)$. A universal property of $M_\gp(N)$ implies that $M \cong V_\gp(N)$. Thus the original classification problem is reduced to classifying all possible weight
$\gl$--modules $N$. This leads to an inductive description of all irreducible weight modules with finite dimensional weight spaces once the 
torsion free modules are known.

The argument above illustrates how the interplay between two properties of parabolic subalgebras is used for the classification problem under consideration.
First we constructed the parabolic subalgebra $\gp$ using the fact that its roots satisfy properties (\ref{P1}) and (\ref{P2}). Then we used the fact that $\gp$ induces a
$\Z$--grading on $\gg$ to construct the module $V_\gp(N)$. 

It is natural to attempt to generalize this scheme beyond the case of finite dimensional reductive Lie algebras.
The first step in such a program is to introduce the notion of parabolic subalgebra. 
For a reductive finite dimensional Lie algebra $\gg$ the definition is intrinsic --- a parabolic
subalgebra is any subalgebra which contains a Borel subalgebra.
Fixing a Cartan subalgebra $\gh$ of $\gg$, the conjugacy theorem implies that every parabolic
subalgebra is conjugate to a parabolic subalgebra $\gp$ containing $\gh$. If $\gh \subset \gp$, 
then $\gp$ is of the form (\ref{parabolic})
for some subset $P \subset \Delta$. The fact that $\gp$ contains a Borel subalgebra is equivalent to 
(\ref{P1}) and the fact that $\gp$ is a subalgebra of $\gg$ is
equivalent to (\ref{P2}). In this way we obtain an alternative definition of a parabolic subalgebra. 
Namely, we can define a parabolic subalgebra containing the fixed
subalgebra $\gh$ as a subalgebra of the form (\ref{parabolic}), where $P$ satisfies (\ref{P1}) and (\ref{P2}). 
A set of roots $P$ with the properties (\ref{P1}) and (\ref{P2}) is called parabolic, \cite{Bo}.

V. Futorny pioneered the study of weight modules over affine Lie algebras.  
In \cite{Fu1} he described explicitly all parabolic sets
of roots for all affine Lie algebras. The corresponding
parabolically induced modules were then studied in 
\cite{Fu3}, \cite{Fu2},  \cite{FKM}, and \cite{FK}. Irreducible weight modules with
finite dimensional weight spaces on which the center acts non trivially
 were classified in \cite{FT}. The
classification of all irreducible weight modules with finite
dimensional weight spaces over affine Lie algebras was recently
completed in \cite{DG}. Not unexpectedly, the first step in the
latter classification is a parabolic induction theorem analogous to the
one for finite dimensional reductive Lie algebras. 

Other cases to which the scheme above was
applied include finite dimensional simple Lie superalgebras,
\cite{DMP}, and affine Lie superalgebras, \cite{FR} and
\cite{Fu4}.

The aim of this note is to use an analog of (\ref{grading})
for introducing the notion of strongly parabolic set of roots. 
We then compare parabolic and strongly parabolic sets of roots for
simple finite dimensional, affine, and toroidal Lie algebras and superalgebras. 
Even though these definitions are not equivalent for some Lie
superalgebras, they are equivalent for all Lie algebras  and for most Lie superalgebras 
under consideration. For finite dimensional Lie algebras and superalgebras
the results were known except for the Lie superalgebra $H(n)$. For
affine Lie algebras the equivalence can be derived from the work
of Futorny, see \cite{Fu2}. Our approach provides a uniform
treatment of both affine Lie algebras and superalgebras. The
results for toroidal Lie algebras and  superalgebras are new. Here
is briefly the content of each of the six sections.

\begin{itemize}
\item[{1.}] Definitions.
\item[{2.}] Kac--Moody Lie algebras and superalgebras.
\item[{3.}] Simple finite dimensional Lie superalgebras.
\item[{4.}] Affine Lie algebras and superalgebras.
\item[{5.}] Toroidal Lie algebras and superalgebras.
\item[{6.}] Conclusion.
\end{itemize}

{\bf Acknowledgment.} I. D. acknowledges the support of FAPESP
and the hospitality and excellent working conditions at the
Institute for Mathematics and Statistics at the University of
S\~ao Paulo, where most of this work was completed.

\section{Definitions} \label{definitions}
Let $V$ be a finite dimensional real vector space and let $\Delta \subset V \setminus \{0\}$.

\begin{definition} $\phantom{xx}$ \label{def}

\begin{itemize}
\item[(i)] A subset $P \subset \Delta$  is called parabolic if it satisfies conditions (\ref{P1}) and (\ref{P2}).

\item[(ii)] A partition $\Delta = \Delta^- \sqcup \Delta^0 \sqcup \Delta^+$ is called  a triangular decomposition of $\Delta$ if there exists a
linear function $\lambda \in V^*$ such that $\Delta^0 = \Delta \cap \Ker \lambda$ and
$\Delta^\pm = \{ \alpha \in \Delta \, | \, \lambda(\alpha) \gtrless 0 \}$.

\item[(iii)] A subset $P \subset \Delta$ is called strongly parabolic if
$P = \Delta$ or $P = P^0 \sqcup \Delta^+$ for some triangular decomposition $\Delta = \Delta^- \sqcup \Delta^0 \sqcup \Delta^+$ and
some strongly parabolic subset $P^0 \subset \Delta^0$ considered in the vector space $\Ker \lambda \subset V$.

\item[(iv)] A strongly parabolic set $P$ is called a principal parabolic set if there exists a triangular decomposition
$\Delta = \Delta^- \sqcup \Delta^0 \sqcup \Delta^+$ such that
$P = \Delta^0 \sqcup \Delta^+$.
\end{itemize}
\end{definition}

The following observation is clear.

\begin{proposition} \label{prop1} If $\Delta = - \Delta$, then
every strongly parabolic subset of $\Delta$ is a parabolic subset of $\Delta$.
\end{proposition}

If $\Delta \neq - \Delta$ the notion of parabolic set does not seems to be the correct generalization of the notion of parabolic
set of roots of finite dimensional reductive Lie algebras. Problems arise already for finite dimensional Lie superalgebras and for
some natural classes of infinite dimensional Lie algebras, see sections \ref{kac} and \ref{conclusion}.
The converse of Proposition \ref{prop1} is not true in general --- see the discussion about $psl(m|m)$ and $H(n)$ in section \ref{kac} below.

For the rest of this paper
$\Delta$ will be the root system of a Lie algebra or superalgebra $\gg$ (with respect to a fixed subalgebra $\gh$) considered in the real vector
space  $V = \R \otimes_\Z Q$, where $Q$ is the abelian group generated by $\Delta$. The case when $\gg$ is a
reductive finite dimensional Lie algebra is classical. The proof of the following statement is standard, see Proposition VI.7.20 in \cite{Bo}.

\begin{proposition} \label{prop2}
Let $\Delta$ be the root system of a reductive finite dimensional Lie algebra $\gg$ and let $P \subset \Delta$. The following are equivalent.
\begin{itemize}
\item[(i)] $P$ is parabolic.
\item[(ii)] $P$ is strongly parabolic.
\item[(iii)] $P$ is principal parabolic.
\end{itemize}
\end{proposition}

Next we show that strongly parabolic sets yield subalgebras and subsuperalgebras with well--behaved parabolic induction functor.
Assume that $\Delta = - \Delta$.
Let $P$ be a strongly parabolic subset of $\Delta$ with corresponding triangular decomposition
$\Delta = \Delta^- \sqcup \Delta^0 \sqcup \Delta^+$
and strongly parabolic subset $P^0$ of $\Delta^0$. Set $\gg^0 = \gh \oplus (\oplus_{\alpha \in \Delta^0} \, \gg^\alpha)$, $\gg^\pm = \oplus_{\alpha \in \Delta^\pm} \gg^\alpha$,
$\gp = \gh \oplus (\oplus_{\alpha \in P} \, \gg^\alpha)$, $\gp^0 = \gp \cap \gg^0$, $\gp^+ = \gp \cap \gg^+$, $\gl = \gh \oplus (\oplus_{\alpha \in P \cap -P} \, \gg^\alpha)$, and
$\gn = \oplus_{\alpha \in P \setminus -P} \, \gg^\alpha$.

\begin{proposition} \label{prop3} Let $\Delta = - \Delta$, let $P \subset \Delta$ be a strongly parabolic subset of $\Delta$,
 and let $\gp$, $\gl$ and $\gn$ be as above. The following hold.

\begin{itemize}
\item[(i)] $\gn$ is an ideal of $\gp$.

\item[(ii)] Every irreducible weight $\gl$--module $N$ is an irreducible weight $\gp$--module with the trivial action of $\gn$. Conversely,
if $N$ is an irreducible weight $\gp$--module, then $\gn$ acts trivially on $N$ and $N$ is an irreducible weight $\gl$--module.

\item[(iii)] If $N$ is an irreducible weight $\gp$--module, then the module $M_\gp(N) = U(\gg) \otimes_{U(\gp)}  N$ admits a unique
maximal proper submodule and, consequently, a unique irreducible quotient $V_\gp(N)$.
\end{itemize}
\end{proposition}

\noindent
{\bf Proof.} If $P = \Delta$, there is nothing to prove. Assume now that $P \subset \Delta$ is proper.

\noindent
(i) Let $\alpha \in P \setminus{-P}$ and $\beta \in P$ with $\alpha + \beta \in \Delta$.  We need to show that $\alpha + \beta \in P \setminus -P$.
If $\alpha \in \Delta^+$ or $\beta \in \Delta^+$, then $\alpha + \beta \in \Delta^+ \subset P \setminus -P$. If both $\alpha, \beta \in \Delta^0$, then
$\alpha + \beta \in P \setminus -P$ by induction.

\noindent
(ii) The fact that every irreducible $\gl$--module has a structure of a $\gp$--module follows from (i).
Let $N$ be an irreducible weight $\gp$--module and let $\mu$ be a weight of $N$. Set
$N' = \oplus_{\lambda(\mu') \geq \lambda(\mu)} N^{\mu'}$ and $N'' = \oplus_{\lambda(\mu'') > \lambda(\mu)} N^{\mu''}$.
Both $N'$ and $N''$ are submodules of $N$, the former is nonzero and the latter is proper. Since $N$ is irreducible, $N' = N$ and $N'' = 0$.
This shows that $\gp^+$ acts trivially on $N$ and $N$ is an irreducible $\gp^0$-module.
An induction argument implies that $\gn$ acts trivially on $N$.

\noindent
(iii) We proceed by induction. Fix a weight $\mu \in \supp N$.
Let $X'$ be the maximal proper submodule of the $\gg^0$--module $U(\gg^0) \otimes_{U(\gp^0)} N$
and let $M'$ be the $\gg$--submodule of $M_\gp(N)$ generated by $X'$. It is clear that $M'$ is a proper submodule of $M_\gp(N)$.
Let $M''$ be the sum of all submodules $Y$ of $M_\gp(N)$ with the property that $\mu' \in \supp Y$ implies $\lambda(\mu') < \lambda(\mu)$.
We leave it to the reader to verify that $M' + M''$ is a proper submodule of $M_\gp(N)$ and that every proper submodule of $M_\gp(N)$
is a submodule of $M' + M''$. \hfill $\square$

\vskip.2in

We conclude this section with the remark that there 
are interesting Lie algebras and superalgebras with infinite dimensional Cartan subalgebras. Their roots generate infinite dimensional
real vector spaces and may require a more general notion of strongly parabolic sets of roots. This is beyond the scope of the present paper.
For a definition of Borel subalgebra of a Lie algebra with infinite dimensional Cartan subalgebra see \cite{DP}.

\section{Kac--Moody Lie algebras and superalgebras} \label{kac}
The purpose of this section is to study certain parabolic sets of roots for a broad class of Lie algebras and superalgebras.
The main result, Proposition \ref{kac-moody}, will provide the relationship between parabolic sets and strongly parabolic sets
of roots for most simple finite dimensional Lie superalgebras with $\Delta = - \Delta$. Proposition \ref{kac-moody} will also allow us
to treat affine Lie algebras and superalgebras in a uniform manner.

Kac--Moody Lie algebras were introduced in the 70's and have been a focus of research since then. The idea behind the
definition is simple --- one considers Lie algebras defined by generators and relations which are encoded in a Cartan matrix.
Defining the analogous Lie superalgebras presents some challenges, most importantly, the existence of odd simple roots
along which there are no reflections in the classical sense. In \cite{K1} Kac introduced contragredient Lie superalgebras
 defined by generators and relations and classified the simple finite dimensional contragredient Lie superalgebras.
  Recently Serganova gave a definition of Kac--Moody Lie superalgebras and jointly with Hoyt classified all
finite--growth contragredient Lie superalgebras. Since the structure theory of Kac--Moody Lie algebras and superalgebras is
beyond the scope of this paper, we will only present the properties that we need as well as the list of simple finite dimensional and
affine Lie superalgebras which are quasisimple regular Kac--Moody Lie superalgebras.
For details on Kac--Moody Lie algebras we refer the reader to \cite{K3},
and on Kac--Moody Lie superalgebras --- to \cite{S2}.

For the rest of this section we assume that $\gg$
is a Kac--Moody Lie algebras or a  quasisimple regular Kac--Moody Lie superalgebra with a fixed Cartan subalgebra $\gh$.
The subalgebra $\gh$ is finite
dimensional and $\gg$ decomposes as in (\ref{eq11}). Since $\Delta = -\Delta$, Proposition \ref{prop1} implies that every
strongly parabolic subset of $\Delta$ is parabolic. Following \cite{S2} we call a subset $\Sigma \subset \Delta$ a {\it base of $\Delta$} if
$\Sigma$ is linearly independent and for every $\alpha \in \Sigma$ there exist elements $X_\alpha \in \gg^\alpha$, $Y_\alpha \in \gg^{-\alpha}$
such that $\{X_\alpha, Y_\alpha\}_{\alpha \in \Sigma} \cup \gh$ generate $\gg$. The elements of $\Sigma$ are called {\it simple roots}.
Every element $\beta \in \Delta$ is an integer combination of elements of $\Sigma$ with nonnegative or nonpositive coefficients only. We denote by $\Delta^+(\Sigma)$ the set of roots which are nonnegative integer combinations of elements of $\Sigma$.
If $\beta$ is a negative root, i.e. $-\beta \in \Delta^+(\Sigma)$, then either $- \beta \in \Sigma$ or
there exists $\alpha \in \Sigma$ such that $\beta + \alpha \in \Delta$ is still a negative root.
If $\alpha \in \Sigma$, then we have one of the following alternatives.
\begin{itemize}
\item[(i)] $\alpha$ is an even root. In this case $k \alpha \in \Delta$ if and only if $k = \pm 1$ and
$X_\alpha, Y_\alpha$ generate a Lie algebra isomorphic to $sl_2$;
\item[(ii)] $\alpha$ is an odd root and $k \alpha \in \Delta$ if and only if $k = \pm 1, \pm 2$. In this case
$X_\alpha, Y_\alpha$ generate a Lie superalgebra isomorphic to $osp(1|2)$;
\item[(iii)] $\alpha$ is an odd root and $k \alpha \in \Delta$ if and only if $k = \pm 1$. In this case
$X_\alpha, Y_\alpha$ generate a Lie superalgebra isomorphic to $sl(1|1)$.
\end{itemize}
For every simple root $\alpha$ we denote by $r_\alpha$ the corresponding reflection. In cases (i) and (ii) $r_\alpha$ is the
usual reflection along the even root $\alpha$ or $2 \alpha$ respectively. In case (iii) $r_\alpha$ is an odd reflection.
Odd reflections were introduced by Penkov and Serganova, \cite{PS}, to compensate for the fact that the Weyl group
of a Lie superalgebra is usually too small. For a detailed treatment of odd reflections, see \cite{S2}.
For every $\alpha \in \Sigma$ the reflection $r_\alpha: \Delta \to \Delta$ is a bijection and $r_\alpha(\alpha) = - \alpha$.
Furthermore, $\Sigma' = r_\alpha(\Sigma)$
is a base of $\Delta$ and $\Delta^+(\Sigma') \setminus \{-\alpha, -2 \alpha\} = \Delta^+(\Sigma) \setminus \{\alpha, 2 \alpha\}$, i.e.
the only positive roots which $r_\alpha$  makes negative are $\alpha$ and, in case (ii), $2 \alpha$.

The following proposition will allow us to establish the relationship between parabolic and strongly parabolic sets of roots for most
simple finite dimensional Lie superalgebras and will be crucial in our treatment of affine Lie algebras and superalgebras.

\begin{proposition} \label{kac-moody} Let $P$ be a parabolic subset of $\Delta$. Assume that there is a base $\Sigma$ of $\Delta$
for which $\Delta^+(\Sigma) \setminus P$ is a finite set. Then

\begin{itemize}
\item[(i)] There exists a base $\Pi$ of $\Delta$ such that $\Delta^+(\Pi) \subset P$;

\item[(ii)] $P$ is a principal parabolic subset of $\Delta$.
\end{itemize}
\end{proposition}

\noindent
{\bf Proof.}   The proof is a variation of the proof of Proposition VI.7.20 in \cite{Bo}.

\noindent
(i) Let $\Pi$ be a base of $\Delta$ such that the cardinality of the set $\Delta^+(\Pi) \setminus P$ is minimal.
If $\Delta^+(\Pi) \not \subset P$, then there is $\alpha \in \Pi$ such that $\alpha \not \in P$. If $\Pi' = r_\alpha(\Pi)$, then one checks easily that
$\Delta^+(\Pi') \setminus P = (\Delta^+(\Pi) \setminus P) \setminus \{\alpha, 2 \alpha \}$ which contradicts the choice of
$\Pi$.

\noindent
(ii) For $\alpha \in \Pi$ set
\begin{equation} \label{lambda}
\lambda(\alpha) = \left\{ \begin{array}{ccl} 1 & {\text { if }} & -\alpha \not \in P\\
0 & {\text { if }} & -\alpha  \in P. \end{array} \right.
\end{equation}
Since $\Pi$ is a basis of $V$, (\ref{lambda}) defines an element $\lambda \in V^*$. Let $\Delta = \Delta^- \sqcup \Delta^0 \sqcup \Delta^+$
be the corresponding triangular decomposition. The fact that $P$ is closed implies that $\Delta^0 \sqcup \Delta^+ \subset P$.
Assume $\Delta^- \cap P \neq \emptyset$ and pick $\beta \in \Delta^- \cap P$ of minimal height, i.e. such that
$\sum_{\alpha \in \Pi} k_\alpha$ is minimal where $\beta = - \sum_{\alpha \in \Pi} k_\alpha \alpha$. By the definition of $\lambda$,
$- \beta \not \in \Pi$. Let $\alpha_0 \in \Pi$ be such that $\beta + \alpha_0 \in \Delta$. Then $\beta + \alpha_0 \in \Delta^- \cap P$
and $\beta + \alpha_0$ has smaller height than $\beta$ which is a contradiction. Thus $\Delta^- \cap P= \emptyset$, i.e.
$P = \Delta^0 \sqcup \Delta^+$. The proof is complete. \hfill $\square$

\section{Simple finite dimensional Lie superalgebras}
 We use the notation from \cite{K1} and \cite{S2}.
The simple finite dimensional Lie superalgebras which are not Lie algebras are: $sl(m|n)$ for $m \neq n$, $psl(m|m)$, $osp(m|2n)$,
$D(2,1;\alpha)$, $F(4)$, $G(3)$, $P(n)$, $psq(n)$, and the Cartan type superalgebras $W(n)$, $S(n)$, $\tilde{S}(n)$, and $H(n)$. For the
restrictions on the parameters $m,n$, and $\alpha$ as well as isomorphisms among the supealgebras listed above
we refer the reader to \cite{K1}. The root systems of the superalgebras $P(n)$, $W(m)$, $S(m)$, and $\tilde{S}(m)$, $n\geq 2$, $m\geq 3$, do not satisfy the condition $\Delta = - \Delta$ and we will not comment on them. Of the rest, $sl(m|n)$ for $m \neq n$,
$osp(m|2n)$, $D(2,1;\alpha)$, $F(4)$, and $G(3)$ are  quasisimple regular Kac--Moody Lie superalgebras and Propostion \ref{kac-moody}
applies. The root system of $psq(n)$ is the same as the root system of the Lie algebra $sl_n$ and hence the notions of parabolic set,
strongly parabolic set, and principal parabolic set for $psq(n)$ and $sl_n$ coincide. In particular every parabolic set of roots
of $psq(n)$ is a principal parabolic set. The remaining cases of $psl(m|m)$ and $H(n)$ are discussed below.

\vskip.2in
\noindent
{\bf The case of $\gg = psl(m|m), m \geq 2$.} The root system is $\Delta = \{ \vep_i - \vep_j, \tau_i - \tau_j \, | \, 1 \leq i \neq j \leq m\} \cup
\{\pm (\vep_i - \tau_j) \, | \, 1 \leq i, j \leq m\} \subset V$, where $V = \Span \{\vep_i, \tau_i \, | \, 1 \leq i \leq j\}$ subject to the relations
 $\vep_1 + \ldots + \vep_m = \tau_1 + \ldots + \tau_m = 0$. If $m = 2$, then
$\Delta$  is the same as the root system $B_2$ of the Lie algebra $so_5$. As a consequence, every parabolic subset of $\Delta$ is a
principal parabolic set.  Note, however, that  $psl(2|2)$ requires special attention as the odd root spaces are
two dimensional. Let now $m>3$. It is easy to see that the set $P = \{ \vep_i - \vep_j, \tau_i - \tau_j \, | \, 1 \leq i < j \leq m\} \cup
\{ \vep_i - \tau_j \, | \, 1 \leq i \neq j \leq m \}$ is a parabolic subset of $\Delta$. It is not, however, a strongly parabolic subset since
\begin{equation} \label{sum}
\begin{array}{l}
(\vep_1 - \vep_2) + 2 (\vep_2 - \vep_3) + \ldots + (m-1) (\vep_{m-1} - \vep_m) + m (\vep_m - \tau_1)  \\ +
(m-1)  (\tau_1 - \tau_2)
+ \ldots + 2 (\tau_{m-2} - \tau_{m-1}) + (\tau_{m-1} - \tau_m) = 0.
\end{array}
\end{equation}
This phenomenon is due to the fact that $psl(m|m)$ is not a Kac--Moody Lie superalgebra. To understand the parabolic subsets of $\Delta$
one needs to consider the  quasisimple regular Kac--Moody Lie superalgebra $gl(m|m)$ corresponding to $psl(m|m)$.
The root system of $gl(m|m)$ is $\Delta' = \{ \vep_i - \vep_j, \tau_i - \tau_j \, | \, 1 \leq i \neq j \leq m\} \cup
\{\pm (\vep_i - \tau_j) \, | \, 1 \leq i, j \leq m\} \subset V'$, where $V' = \Span \{\vep_i, \tau_i \, | \, 1 \leq i \leq j\}$ without any
relations among the generators. Proposition \ref{kac-moody} applies to $\Delta'$, i.e.
every parabolic subset of $\Delta'$ is a principal parabolic set. It is clear now that the parabolic
subsets of $\Delta$ are exactly the images of the principal parabolic subsets of $\Delta'$ under  the obvious surjection $V' \to V$.
Note that (\ref{sum}) implies that even for the supealgebra $sl(m|m)$ not every parabolic subset of roots is strongly parabolic.

\vskip.2in
\noindent
{\bf The case of $\gg = H(n), n \geq 5$.} The root system is $\Delta = \{ \sum_{i=1}^l k_i \vep_i \, | \, k_i \in \{0, \pm 1\}\}$, where $l = \lfloor n/2 \rfloor$ is
the integer part of $n/2$ and $V = \Span \{\vep_1, \ldots, \vep_l\}$. The root system of $H(2l)$ is the same as the root system of $H(2l+1)$ and the
root system of $H(5)$ is the same as the root system of the Lie algebra $so_5$. We do not know whether every parabolic subset of $\Delta$ is
strongly parabolic for $6 \leq n \leq 9$. For $n \geq 10$, however, there are parabolic subsets of $\Delta$ which are not strongly parabolic.

Here is an example for $n =10$, i.e. $l=5$. If $\alpha = \sum_{i=1}^{5} k_i \vep_i$ is a root, set 
$\hgt(\alpha) = \sum_{i=1}^{5} k_i$ and define
\begin{align*}
P_0  = & \{ \alpha \in \Delta  \, | \, \hgt(\alpha) = 0,  k_1 \geq 0, k_1 + k_2 \geq 0, k_1 + k_2 + k_3 \geq 0, k_1 + k_2 + k_3 + k_4 \geq 0\} \\
 \cup & \{-\vep_1+ \vep_2 + \vep_3 - \vep_4, -\vep_1+ \vep_2 + \vep_3 - \vep_5, -\vep_1+ \vep_2 + \vep_4 - \vep_5, \\ & -\vep_1+ \vep_3 + \vep_4 - \vep_5,
-\vep_2+ \vep_3 + \vep_4 - \vep_5\} {\text { and }}
\end{align*}
$$
P = P_0 \cup \{\alpha \in \Delta \, | \, \hgt(\alpha) >0\}.
$$
A direct verification shows that $P$ is a parabolic subset of $\Delta$. Assume that $P$ is a strongly parabolic subset of $\Delta$ with corresponding
linear function $\lambda \in V^*$. Then $\lambda$ vanishes on each of the roots $\{-\vep_1+ \vep_2 + \vep_3 - \vep_4, -\vep_1+ \vep_2 +
\vep_3 - \vep_5, -\vep_1+ \vep_2 + \vep_4 - \vep_5,  -\vep_1+ \vep_3 + \vep_4 - \vep_5,
-\vep_2+ \vep_3 + \vep_4 - \vep_5\}$, which implies that $\lambda(\vep_1) = \lambda(\vep_2) = \ldots = \lambda(\vep_5) >0$. Hence
$\Delta^0 = \{ \alpha \in \Delta \, | \, \hgt(\alpha) = 0 \}$ and $P_0 = P \cap \Delta^0$ is a strongly parabolic subset of $\Delta^0$.
The vector space $V^0$ is spanned by $\{\vep_1 - \vep_2, \vep_2 - \vep_3, \vep_3 - \vep_4, \vep_4 - \vep_5\}$. Let $\lambda^0 \in (V^0)^*$
be the linear function corresponding to the strongly parabolic subset $P_0$ of $\Delta^0$.  Then $\lambda^0$ vanishes on
each of the roots $\{-\vep_1+ \vep_2 + \vep_3 - \vep_4, -\vep_1+ \vep_2 +
\vep_3 - \vep_5, -\vep_1+ \vep_2 + \vep_4 - \vep_5,  -\vep_1+ \vep_3 + \vep_4 - \vep_5,
-\vep_2+ \vep_3 + \vep_4 - \vep_5\}$, which, as the reader will verify, implies that $\lambda^0 = 0$. This contradicts the fact that $P_0$ is a
proper subset of $\Delta^0$. Hence $P$ is not a strongly parabolic subset of $\Delta$.

The example above can be extended to an example of a parabolic subset of $\Delta$ which is not strongly parabolic for any $n > 10$ in
the following way. Consider the natural projection $\pi: \Delta(n) \to \Delta(10) \cup\{0\}$, where $\Delta(n)$ denotes the roots of $H(n)$.
The set $\pi^{-1}(P \cup \{0\})$, where $P \subset \Delta(10)$ is as above, is a parabolic subset of $\Delta(n)$ which is not strongly parabolic.

\section{Affine Lie algebras and superalgebras} \label{affine}

Let $\gs$ be a finite dimensional Lie algebra or superalgebra not isomorphic to $psl(m|m)$ 
with an invariant bilinear form $\kappa$. The affine Lie algebra
or superalgebra $\gs^{(1)}$ is defined as
$$
\gs^{(1)} = \gs \otimes \C[t, t^{-1}] \oplus \C D \oplus \C K
$$
with commutation relations
$$
[x \otimes t^k, y \otimes t^l] = [x,y] \otimes t^{k+l} + \delta_{k, -l} k \kappa(x,y) K, \quad
[D, x \otimes t^k] = k x \otimes t^k, \quad [K, \gs^{(1)}] = 0,
$$
where $\delta_{k,-l}$ is Kronecker's delta function.
If $\gs$ is a simple finite dimensional Lie algebra, then $\gs^{(1)}$ is a Kac--Moody Lie algebra. If $\gs$ is
one of the superalgebras $sl(m|n)$ for $m \neq n$, $osp(m|2n)$, $D(2,1;\alpha)$, $F(4)$, or $G(3)$, then
$\gs^{(1)}$ is a quasisimple regular Kac--Moody superalgebra. Finally, we define $psl(m|m)^{(1)}$ for $m \geq 3$
as
$$
psl(m|m)^{(1)} = psl(m|m) \otimes \C[t, t^{-1}] \oplus \C D \oplus \C K \oplus \C D' \oplus \C K'
$$
with commutation relations
$$
[x \otimes t^k, y \otimes t^l] = [x,y] \otimes t^{k+l} + \delta_{k, -l} k \kappa(x,y) K + \delta_{k,-l} \tr([x,y]) K',
$$
$$
[D, x \otimes t^k] = k x \otimes t^k, \quad [D', x \otimes t^k] = (1 - (-1)^{p(x)}) x \otimes t^k,
$$
$$
[K, psl(m|m)^{(1)}] = 0,
\quad [K', psl(m|m)^{(1)}] = 0,
$$
where $p(x)$ denotes the parity of $x$. The Lie superalgebra $psl(m|m)^{(1)}$ is a quasisimple regular Kac--Moody
Lie superalgebra. (Note that $psl(m|m)^{(1)}$ is a subalgebra of $gl(m|m)^{(1)}$ with the same root system. In order to
have a uniform treatment of affine Lie algebras and Lie superalgebras we will consider $psl(m|m)^{(1)}$ instead of
the more natural superalgebra $gl(m|m)^{(1)}$.)
For the rest of this section $\gs^{(1)}$ will denote one the Lie algebras or superalgebras above.

Denote by $R$ the roots of $\gs$ and by $W$ the real vector space spanned by $R$. The root system $\Delta$ of
$\gs^{(1)}$ is given by
$$
\Delta =  (R + \Z \delta) \cup (\Z \setminus \{0\}) \delta,
$$
where $\delta$ is one of the two indivisible imaginary roots of $\gs^{(1)}$. The real vector space spanned by $\Delta$
is $V = W \oplus \R \delta$.

If $\theta$ is an automorphism of $\gs$ of order $p$ which preserves $\kappa$ and $\vep$ is a p$^{th}$ primitive root if $1$,
we extend $\theta$ to an automorphism of $\gs^{(1)}$ by setting
$$
\theta(x \otimes t^k) = \vep^k \theta(x) \otimes t^k,
$$
$$
\theta(D) = D, \quad \theta(K) = K, \quad \theta(D') = D', \quad \theta(K') = K'.
$$
We denote the fixed points of $\theta$ by $\gs^{(p)}$. The so obtained Lie algebra (or superalgebra) is a twisted affine Lie algebra
(or superalgebra). (This notation is somewhat ambiguous since the isomorphism class of
a twisted affine Lie superalgebra does not determine the order of $\theta$. Since we limit our considerations to
the list of twisted affine Lie algebras and superalgebras, we will not run into problems.) The twisted affine Lie algebras are $A_n^{(2)}, D_n^{(2)}$,
$E_6^{(2)}$, and $D_4^{(3)}$ and they are Kac--Moody Lie algebras. Among the twisted affine Lie superalgebras, the quasisimple regular Kac--Moody Lie superalgebras are:
$sl(m|n)^{(2)}$ for even $mn$, $psl(m|m)^{(2)}$ for even $m$, $sl(m|n)^{(4)}$ for odd $mn$, $psl(m|m)^{(4)}$ for odd $m$,
$osp(2m|2n)^{(2)}$, and $q(n)^{(2)}$. It is clear from the
definition of $R$ that $\Delta \subset (R + \Z \delta) \cup (\Z \setminus \{0\}) \delta$ but, unlike the untwisted case, the inclusion
is proper.

For the twisted affine Lie algebras and superalgebras we consider again the roots $\Delta \subset V$, write
$V = W \oplus \R \delta$ and define $R$ as the image of $\Delta \setminus \Z \delta$ under the natural projection $V \to W$. The set
$R \subset W$ is again a root system, which may no longer be reduced. For example, for $A_{2l}^{(2)}$, $R$ is the nonreduced
root system $BC_l$. It is not difficult, however, to check that every parabolic subset of $R$ is a principal parabolic set.

For the rest of the section $\gg$ will denote one of the algebras or superalgebras $\gs^{(p)}$ above.

\begin{theorem} \label{prop_affine} If  $P \subset \Delta$ is a parabolic set of roots, then $P$ is strongly parabolic.
Furthermore, if $P$ is a proper subset of $\Delta$, one of the three mutually exclusive alternatives holds.

\begin{itemize}
\item[(i)] $P$ is principal and $\delta \not \in \Delta^0$;

\item[(ii)] $P$ is principal and $(\Z \setminus \{0\}) \delta \subset P$;

\item[(iii)] $P$ is not principal.
\end{itemize}
\end{theorem}

\noindent
{\bf Proof.} If $P = \Delta$ then $P$ is strongly parabolic by definition. Assume that $P$ is a proper subset of $\Delta$.
We first prove that $P$ is strongly parabolic by considering two cases.

\noindent
{\bf Case 1.} There exists $\alpha \in R$ such that $( \alpha + \Z \delta ) \cap \Delta \subset P$. Set
$$
S = \{ \alpha \in R \, | \, (\alpha + \Z \delta) \cap P \neq \emptyset \}.
$$
Then $S$ is a parabolic subset of $R$. This is obvious in the untwisted case and requires an elementary check in
the twisted case. Furthermore, $S$ is a proper subset of $R$. Using the results from the previous section we conclude
that $S$ is a principal parabolic subset of $R$. Let $\lambda \in W^*$ be a corresponding linear function.
Extend $\lambda$ to a linear function $\Lambda \in V^*$ by setting $\Lambda(\delta) = 0$ and consider the triangular
decomposition of $\Delta$ corresponding to $\Lambda$. From the definition of $S$ it follows that $\Delta^- \cap P = \emptyset$,
which also implies that $\Delta^+ \subset P$. Finally, $P \cap \Delta^0$ is a parabolic subset of $\Delta^0$ and a
simple inductive argument completes the proof that $P$ is strongly parabolic.

\noindent
{\bf Case 2.} For every $\alpha \in R$, $(\alpha + \Z \delta) \cap P \neq \emptyset$. Since either $\delta$ or $-\delta$ belongs to
$P$, we may assume that $\delta \in P$. This implies that for every $\alpha \in R$, $\alpha + n \delta \in P$ for large enough $n$ and, consequently,
 $\Delta^+ (\Sigma) \setminus P$ is a finite set, where $\Sigma$ is the standard base of $\Delta$.
Applying Proposition \ref{kac-moody} we obtain that $P$ is a principal parabolic set and $\delta \not \in \Delta^0$.

To complete the proof we notice that in Case 1 we have two alternatives for the parabolic subset $P \cap \Delta^0$ of $\Delta^0$.
Either $P \cap \Delta^0 = \Delta^0$ and then $P$ is principal and $(\Z \setminus \{0\}) \delta \subset P$; or
$P \cap \Delta^0$ is a proper subset of $\Delta^0$ and it that case $P$ is not a principal parabolic subset of $\Delta$. However, it is
not difficult to see that in the latter case $P \cap \Delta^0$ is a principal parabolic subset of $\Delta^0$,
i.e. we can think of the parabolic subset from (iii) as "two--step" parabolic subsets.
\hfill $\square$

\begin{definition} \label{def2}
The parabolic sets from cases (i), (ii), and (iii) of Proposition \ref{prop_affine} are called standard, imaginary, and mixed type parabolic sets respectively.
\end{definition}

\section{Toroidal Lie algebras and superalgebras}
Let $\gs$ be a simple finite dimensional Lie algebra or a finite dimensional quasisimple regular Kac--Moody Lie superalgebra, e.g.
$sl(m|n)$ for $m \neq n$, $gl(m|m)$, $osp(m|2n)$, $D(2,1;\alpha)$, $F(4)$, or $G(3)$. For $n \geq 2$ we
define the toroidal Lie algebra or superalgebra $\cT_n(\gs)$ as
\begin{equation} \label{toroidal}
\cT_n(\gs) = \gs \otimes \C[t_1^{\pm 1}, \ldots, t_n^{\pm 1}] \oplus \C D_1 \oplus \ldots \oplus \C D_n \oplus \C K_1 \oplus \ldots \oplus \C K_n
\end{equation}
with commutation relations
\begin{align*}
 [x \otimes t_1^{k_1} \ldots t_n^{k_n}, y \otimes t_1^{l_1} \ldots t_n^{l_n}] & = [x,y] \otimes t_1^{k_1 + l_1} \ldots t_n^{k_n + l_n} \\
 & + \delta_{k_1, -l_1} \ldots \delta_{k_n, - l_n} \kappa(x,y) (k_1 K_1 + \ldots + k_n K_n),
\end{align*}
$$
[D_i, x \otimes t_1^{k_1} \ldots t_n^{k_n}] = k_i x \otimes t_1^{k_1} \ldots t_n^{k_n},  [K_i , \cT_n(\gs)] = 0.
$$
It is easy to check that $\cT_n(\gs)$ admits a root decomposition (\ref{eq11}) with respect to the subalgebra
$\gh \oplus \C D_1 \oplus \ldots \oplus \C D_n \oplus \C K_1 \oplus \ldots \oplus \C K_n$, where $\gh$ is a Cartan subalgebra of $\gs$.
If we denote the roots of $\gs$ by $R$, then the roots of $\cT_n(\gs)$ are given by
\begin{equation} \label{root_toroidal}
\Delta = (R + \Gamma) \cup (\Gamma \setminus \{0\} ),
\end{equation}
where $\Gamma = \Z \delta_1 + \ldots + \Z \delta_n \cong \Z^n$ is the free abelian group with generators $\delta_1, \ldots, \delta_n$,
the elements of $\gh^*$ dual to $D_1, \ldots, D_n$. Furthermore, $V = W \oplus U$, where $W$ is the
real vector space generated by $R$ and $U = \R D_1 \oplus \ldots \oplus \R D_n$.

\begin{theorem} \label{prop_toroidal}
Every parabolic subset of $\Delta$ is strongly parabolic.
\end{theorem}

\noindent
{\bf Proof.} If $P = \Delta$, there is nothing to prove. Assume that $P \subset \Delta$ is proper.
 Following the idea of the proof of Theorem
\ref{prop_affine} we consider two cases.

\noindent
{\bf Case 1.} There exists $\alpha \in R$ such that $\alpha + \Gamma \subset P$. Set
$$
S = \{ \alpha \in R \, | \, (\alpha + \Gamma) \cap P \neq \emptyset \}.
$$
Again $S$ is a proper parabolic subset of $R$ and there exists a triangular decomposition $R = R^- \sqcup R^0 \sqcup R^+$ with
a corresponding linear function $\lambda \in W^*$ for which $S = R^0 \sqcup R^+$. Extend $\lambda$ to a linear function $\Lambda \in V^*$
by setting $\Lambda(\delta_i) = 0$. Let $\Delta = \Delta^- \sqcup \Delta^0 \sqcup \Delta^+$ be the triangular decomposition corresponding to
$\Lambda$. It is clear now that $P = (P \cap \Delta^0) \sqcup \Delta^+$ and $ (P \cap \Delta^0) \subset \Delta^0$ is a parabolic subset, which
completes the proof that $P$ is a strongly parabolic subset of $\Delta$.

\noindent
{\bf Case 2.} For every $\alpha \in R$, $(\alpha + \Gamma) \cap P \neq \emptyset$. Set $T = P \cap \Gamma$. It is clear that $T$ is a proper
parabolic subset of $\Gamma$. A standard separation theorem from analysis implies that there exists a unique linear function $\lambda \in U^*$
such that $T \subset \Gamma^0 \sqcup \Gamma^+$, where $\Gamma = \Gamma^- \sqcup \Gamma^0 \sqcup \Gamma^+$ is the triangular
decomposition corresponding to $\lambda$.
The image of $\lambda$, $\Img \lambda \subset \R$ is either dense or equals $\Z \eta$ for some
$\eta \neq 0$.

If $\Img \lambda$ is dense in $\R$, we extend $\lambda$ to a function $\Lambda \in V^*$ by setting
\begin{equation} \label{eq21}
\Lambda(\alpha) = - \inf \{ \lambda(\gamma) \, | \, \gamma \in \Gamma {\text { with }} \alpha + \gamma \in  P\}.
\end{equation}
It is an easy computation to verify that $\Lambda$ is well--defined and $P \subset \Delta^0 \sqcup \Delta^+$ for the
corresponding triangular decomposition $\Delta = \Delta^- \sqcup \Delta^0 \sqcup \Delta^+$.

Finally, assume that $\Img \lambda = \Z \eta$. We  extend
$\lambda: U \to \R$ to a linear transformation $\theta: V = W
\oplus U \to W \oplus \R \eta$ by setting $\theta_{| W}$ to be the
identity of $W$. The image of $\Delta$ under $\theta$ equals $(R +
\Z \eta) \cup \Z \eta$, i.e. it  is the root system of $\gs^{(1)}$
with the vector $\{0\}$ added to it. Furthermore, since $\theta(P)
\setminus \{0\}$ is a parabolic subset of $\Delta' = (R + \Z \eta)
\cup (\Z \setminus \{0\}) \eta$ we can apply Theorem
\ref{prop_affine} and obtain a linear function $\mu \in (W \oplus
\R \eta)^*$ with corresponding triangular decomposition $\Delta' =
(\Delta')^- \sqcup (\Delta')^0 \sqcup (\Delta')^+$ such that
$\theta(P) \subset (\Delta')^0 \sqcup (\Delta')^+$. Assume that
$\mu(\eta) = 1$ and  extend $\lambda \in U^*$ to a linear function
$\Lambda \in V^*$ by setting $\Lambda(\alpha) = \mu(\alpha)$ for
$\alpha \in R$. It is easy to check that the triangular
decomposition $\Delta = \Delta^- \sqcup \Delta^0 \sqcup \Delta^+$
corresponding to $\Lambda$ satisfies $P \subset \Delta^0 \sqcup
\Delta^+$ and $P \cap \Delta^0$ is a parabolic subset of
$\Delta^0$. \hfill $\square$

\vskip.2in 

\noindent There are some variations in the literature
about the definition of  toroidal Lie algebras. However, all
$\Z^n$--graded central extensions of $\gs \otimes \C[t_1^{\pm1},
\ldots, t_n^{\pm 1}]$ have the same roots system as $\cT_n(\gs)$
and hence Theorem \ref{prop_toroidal} applies to them too.

\section{Conclusion} \label{conclusion}
We established an equivalence of two commonly used combinatorial definitions for the set of roots of parabolic subalgebras
for several classes of Lie algebras and superalgebras. We excluded from our considerations some important Lie superalgebras.
The first class of superalgebras we did not consider are the simple finite dimensional superalgebras whose root systems are not
symmetric, i.e. $\Delta \neq - \Delta$. In these cases the notion of parabolic set does not yield some very natural candidates and
it needs to be generalized. A possible generalization of this notion is the following.

\begin{definition} \label{symmetric}
A subset $P$ of $\Delta$ is parabolic, if $P = P' \cap \Delta$ for some parabolic subset $P'$ of the (symmetric) set $\Delta \cup - \Delta$.
\end{definition}

It is an interesting problem to compare the parabolic sets according to Definition \ref{symmetric} with the strongly
parabolic and principal parabolic sets for the Cartan type Lie superalgebras $W(n)$, $S(n)$, $\widetilde{S} (n)$, and for  $P(n)$. We also note
that there are other definitions of parabolic subalgebras in the literature, cf. \cite{IO} where a $\Z$--grading is used instead of subsets of roots.

We did not consider affine Lie superalgebras which are not quasisimple regular Kac--Moody Lie superalgebras according to \cite{S2}.
This left out superalgebras like $psq(n)^{(1)}$ whose roots are the same as the roots of $sl(n)^{(1)}$ and hence we can use the results
for $sl(n)^{(1)}$. More importantly, we did not consider all twisted affine Lie superalgebras. Indeed, there are more twisted affine Lie superalgebras,
\cite{S1}, \cite{GP}, which however do not have a center. Nevertheless, their representation theory is interesting and may require results about
the relationship between parabolic sets and strongly parabolic sets.

Further classes of Lie algebras and superalgebras for which parabolic subsets are of interest include infinite dimensional
Cartan type Lie algebras and superalgebras, twisted toroidal Lie algebras and superalgebras, Extended Affine Lie Algebras, etc.

\vskip.2in

\noindent
{\bf Note added in proof.} After the paper was accepted for publication we learned about the treatment of parabolic sets in \cite{LN}.
In Chapter 10 of \cite{LN} O. Loos and E. Neher study subsets $P$ of symmetric sets $R$ in real vector spaces. 
Their definition of parabolic sets is stronger
than ours but the two coincide when $R$ satisfies the so called partial sum property. 
Furthermore, the authors introduce  (without a name) a class of parabolic subsets $P$ of $R$
defined via linear maps into partially ordered vector spaces. It is not difficult to verify that 
the class introduced by Loos and Neher coincides with the class of strongly parabolic sets.
Furthermore, Proposition 10.17 in \cite{LN} establishes that the class of
parabolic subsets introduced via linear maps into partially ordered vector spaces exhausts all parabolic sets in the sense of \cite{LN}.
This observation together with the fact that the root systems considered in Theorem \ref{prop_affine} and Theorem \ref{prop_toroidal} 
satisfy the partial sum property
provide  an alternative and unified proof of  the fact that in these cases every parabolic set is strongly parabolic.
 We thank E. Neher for turning our attention to the book \cite{LN} and for suggesting that the class of parabolic sets
introduced there via linear maps into partially ordered vector spaces and the class of strongly parabolic sets coincide.

\vskip.2in

\centerline{\begin{tabular}{lll}
I.D.: & V.F.: & D.G.:\\
Dept. of Math. and Stats.& IME-USP & Dept. of Math.\\
Queen's University& Caixa Postal 66281 & University of Texas at Arlington\\
Kingston, K7L 3N6 & CEP 05315-970, S\~ao Paulo & Arlington, TX 76019\\
Canada & Brazil & USA\\
{\tt dimitrov@mast.queensu.ca} & {\tt futorny@ime.usp.br} & {\tt grandim@uta.edu}
\end{tabular}}

\end{document}